\documentclass[12pt]{article}
\usepackage{amssymb, amsthm, amsmath,eucal,
pstricks,pst-plot,pst-node,eucal}
\usepackage[dvips]{graphics}

\newpsobject{showgrid}{psgrid}{subgriddiv=1,griddots=5,gridlabels=0pt}

\newcommand{\Z}{\mathbb{Z}}
\newcommand{\R}{\mathbb{R}}
\newcommand{\C}{\mathbb{C}}
\renewcommand{\P}{\mathbb{P}}

\newcommand{\bK}{\mathsf{K}}
\newcommand{\cA}{\mathcal{A}}
\newcommand{\cL}{\mathcal{L}}
\newcommand{\cO}{\mathcal{O}}
\renewcommand{\b}{\mathsf{b}}
\newcommand{\w}{\mathsf{w}}

\newcommand{\tC}{\widetilde{C}}
\newcommand{\tP}{\widetilde{P}}
\newcommand{\tbK}{\widetilde{\bK}}

\newcommand{\Rc}{R^\vee}

\DeclareMathOperator{\Coker}{Coker}
\DeclareMathOperator{\Ker}{Ker}
\DeclareMathOperator{\weight}{weight}
\DeclareMathOperator{\Vol}{Vol}
\DeclareMathOperator{\Img}{Img}
\DeclareMathOperator{\tr}{tr}

\numberwithin{equation}{section}

\newtheorem{Theorem}{Theorem}
\newtheorem{Lemma}{Lemma}
\newtheorem{Corollary}[Lemma]{Corollary}
\newtheorem{Proposition}[Lemma]{Proposition}

\begin{document}
\title{Planar dimers and Harnack curves}
\author{Richard Kenyon and Andrei Okounkov}
\date{October 2003}
\maketitle

\section{Introduction}

\subsection{Summary of results}

In this paper we study  
the connection between dimers and Harnack 
curves discovered in \cite{KOS}. To any 
periodic edge-weighted planar bipartite
graph $\Gamma$ one associates its \emph{spectral 
curve} $P(z,w)=0$. The real polynomial $P$
defining the spectral curve arises as the 
characteristic polynomial of the Kasteleyn 
operator in the dimer model. 

It was shown in \cite{KOS} that for real 
positive edge weights on $\Gamma$ the curves
thus obtained are real curves
of a very special kind, namely they
are \emph{Harnack curves}. Harnack 
curves are, in some sense, the
best possible real plane
curves. They were studied both 
classically and recently, 
see \cite{Mi,MR} and references
therein. Here we prove that every Harnack 
curve arises as a spectral curve of some
dimer model.  This gives a 
parameterization of the set of Harnack 
curves which, in spirit, is similar to 
the classical parametrization of totally 
positive matrices. It may be compared
with the result of Vinnikov \cite{Vin}
who gave a similar description of real 
plane curves with a maximally nested 
set of ovals (which form a class of 
curves in some sense opposite to Harnack 
curves). 

We also prove that modulo the natural $(\R^\times)^2$-action
the set of degree $d$ Harnack curves 
in $\P^2(\R)$ is
diffeomorphic to the closed octant 
$\R_{\ge 0}^{(d+4)(d-1)/2}$.
In fact, the areas of the amoeba holes and 
the distances between the amoeba tentacles give
these global coordinates. 

The Kasteleyn operator of the dimer model is
an example of a periodic finite-difference
operator (weighted adjacency matrix of a periodic graph). 
The spectral theory of such operators
is much developed, stimulated, in particular,
by their connections with integrable systems, 
see, for example, \cite{BBT,Dy,GKT} for an 
introduction.  The particular type of 
operators  we consider was studied by A.~Oblomkov \cite{Obl}
and in a series of paper by I.~Dynnikov and S.~Novikov 
\cite{DN1,DN2}. The {\it spectral data}
associated to a periodic finite-difference
operator is, typically, an algebraic 
variety (a curve $C$, in our case) together 
with a line bundle on it, that is, a together 
with a point of the Jacobian $J(C)$. 
With complex coefficients, it is usually easy 
to see (see e.g.\ \cite{CT}) that the ``spectral 
transform'', from the operator to its spectral data,
is dominating, but not surjective. 
Reality issues are typically more subtle. 
{}From the probabilistic origin of our 
spectral  problem, it is natural to consider
real and positive coefficients (edge weights). 
This adds a new aspect to the problem. 

The Harnack curves of genus zero 
play a special role. We characterize
them as the spectral curves of
\emph{isoradial} 
dimers studied in \cite{Ke2}, see Section \ref{sISO},
and also as 
those Harnack curves that minimize 
the volume under their Ronkin function
(with given boundary conditions), 
see Proposition \ref{volmin}. 
Translated into the language of probability, 
this means that isoradial dimer weights
maximize the partition function with 
given boundary conditions.
\bigskip

\subsection{Acknowledgments}

This research was started at the Institut Henri Poincar\'e and 
completed while R.~K.\ was visiting
Princeton University. A.~O.\ was partially supported by
DMS-0096246 and a fellowship from Packard foundation.
 
\section{Background}
\subsection{Kasteleyn operator and its spectral curve}

\subsubsection{}

Let $\Gamma$ be a periodic planar bipartite
graph. We can assume that $\Gamma$ is embedded in the 
plane $\R^2$ in such a way that translations by 
the standard lattice $\Z^2$ preserve $\Gamma$,
including the partition of its vertices into 
black and white ones.  Throughout this paper,
we assume that the 
quotient $\Gamma/\Z^2$ is finite. Let the 
edges of $\Gamma$ be weighted in a
 $\Z^2$--invariant way. 

A \emph{perfect matching} (also known as a dimer 
configuration) of a graph is a collection of 
edges that cover every vertex exactly once. 
The weight of such a matching is defined to 
be the product of the edge weights. 
Kasteleyn in \cite{Kast} 
computed the weighted number of perfect 
matchings of any finite bipartite planar graph $\Gamma_0$ 
using the determinant of what is now called  
the \emph{Kasteleyn 
operator}
\begin{equation}
  \label{Kas}
\bK: \C^{\textup{black vertices}} \to \C^{\textup{white vertices}}\,,
\end{equation}
which is the weighted adjacency matrix of $\Gamma_0$
sign-twisted in a certain way, see \cite{Ke} for 
an introduction. 

The study of perfect matchings (also known as the dimer
model) on infinite periodic planar graphs $\Gamma$
stimulates the study of the corresponding Kasteleyn 
operators, see for example \cite{KOS} for further
information. In this paper, we will focus on 
the case when $\Gamma$ is the hexagonal lattice,
in which case $\bK$ is simply the weighted adjacency matrix. 
This case is, in fact, universal as will be explained
in Section \ref{hexsuffices}.

\subsubsection{}

The operator $\bK$ commutes with the translation 
action of $\Z^2$
and, in particular, it preserves the 
$\Z^2$--eigenspaces.  These eigenspaces
are indexed by characters of $\Z^2$, that is, 
by a pair of Bloch-Floquet 
multipliers $(z,w)\in (\C^*)^2$. They 
are finite--dimensional with a 
distinguished basis $\{\delta_v\}$
consisting of functions supported on a single 
$\Z^2$--orbit and taking value $1$ on a vertex $v$
inside a fixed fundamental domain. Let $\bK(z,w)$ be the 
matrix of $\bK$ in the basis  $\{\delta_v\}$ and acting on the 
$(z,w)$-eigenspace of $\Z^2$ and set, 
by definition, 
\begin{equation}
  \label{defP}
  P(z,w) = \det \bK(z,w) \,.
\end{equation}
Different choices of the fundamental domain 
lead to polynomials that differ by a 
factor of the form $z^i w^j$. 
In particular, the zero locus of $P$
\begin{equation}
  \label{spC}
  \{P(z,w)=0\} \subset (\C^*)^2
\end{equation}
is defined canonically and 
is called the \emph{spectral curve} of
the dimer problem. 

The spectral curve 
remains the same if we multiply 
the weights of all edges incident to 
a given vertex by the same number. 
This is a \emph{gauge
transformations} of the dimer problem;
it does not change probabilities of
configurations. 

The map from the weights modulo gauge
to the corresponding spectral curves
is our main object of study in this  
paper.

\subsubsection{}

\begin{figure}[hbtp]\psset{unit=0.5 cm}
  \begin{center}
    \begin{pspicture}(-2,0)(12,10)
    \rput(5,5){\scalebox{0.64}{\includegraphics{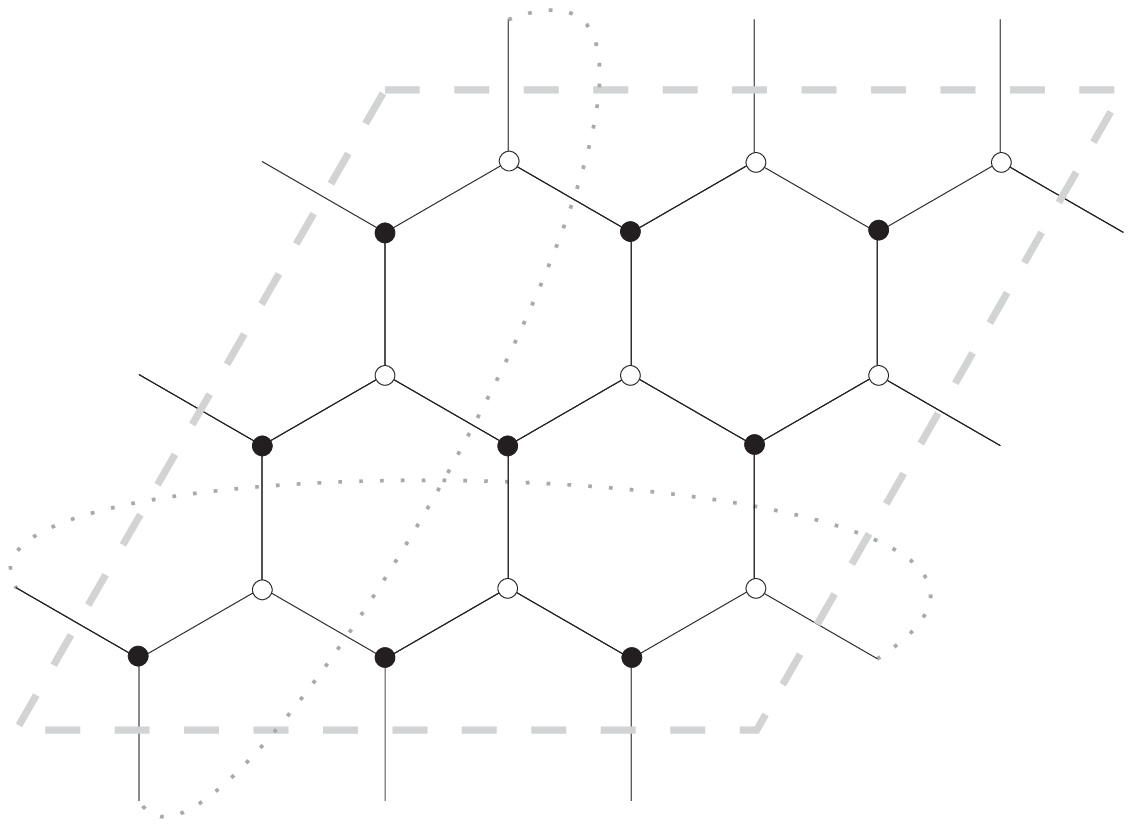}}} 
    \rput[rb](4,9.4){$w$}
    \rput[rb](7.2,9.4){$w$}
    \rput[rb](10.4,9.4){$w$}
    \rput[rt](8.8,1.8){$z$}
    \rput[rt](10.4,4.6){$z$}
    \rput[rt](12,7.3){$z$}
    \rput[rt](5.7,6){$v$}
    \rput[lB](6.59,5.1){$a_v$}
    \rput[l](6,6.5){$b_v$}
    \rput[lt](5,5){$c_v$}
    \end{pspicture}
    \caption{The operator $\bK(z,w)$}
    \label{f1}
  \end{center}
\end{figure}

For example, for the hexagonal lattice with 
$d\times d$-fundamental domain, $P(z,w)$ is
the determinant of the $d^2\times d^2$ 
matrix the  construction of which 
is illustrated in Figure \ref{f1}. 
We denote the weights of the 3 edges incident
to a white vertex $v$ by $a_v$, $b_v$, and $c_v$. 
The dashed line in Figure \ref{f1} is the 
boundary of the fundamental domain. An edge 
crossing it is weighted by an extra factor of $z$ or
 $w$. 

It is clear that in this case $P(z,w)$ is
polynomial in $z$ and $w$ of degree at most
$d$ in each variable. In fact $P$ has total degree $d$:
this can be seen by splitting the fundamental
domain in half into two equilateral triangular arrays of vertices;
the upper right triangle has $d$ more white vertices than
black vertices so exactly $d$ edges connect this
part with the rest of the graph. 

\subsubsection{}

\label{hexsuffices}
The case of the hexagonal lattice is
 universal in the following sense.
By choosing $d$ large enough and setting
some edge weights to zero, one can produce
a graph which is equivalent to an 
arbitrary planar periodic bipartite
graph, in the following sense.

We consider two graphs to be equivalent 
if one can be obtained from another by 
a sequence of moves of the following type:
remove a $1$--valent vertex and its neighbor, 
or remove a $2$--valent vertex, gluing its neighbors into
a single vertex and redistributing
the edge weights accordingly. See \cite{Propp}. It is 
easy to see that such transformation 
induce a weight-preserving bijection of
sets of dimer configurations. 

Because of this universality we will focus
in this paper on the case of hexagonal lattice with 
$d\times d$-fundamental domain. In this case, 
the spectral curve is a degree $d$ curve in 
the projective plane $\P^2$.

\subsection{Harnack curves}

\subsubsection{Topological configuration of ovals}\label{sdefHar}

A real algebraic curve $C(\R)\subset \P^2(\R)$ 
is called an $M$-curve if it has the maximal 
possible number of connected components, namely $1+(d-1)(d-2)/2$,
where $d$ is the degree of $C$. For brevity, we will call 
all these components \emph{ovals}, even though this is at odds
with the classical distinction between separating and nonseparating
components of $C(\R)$ (which will play no role in this paper).
We will also treat isolated real points as degenerate ovals. 
Ovals not intersecting the coordinate axes will be called
\emph{compact ovals}. 

\begin{figure}[!hbtp]
  \begin{center}
    \scalebox{0.54}{\includegraphics{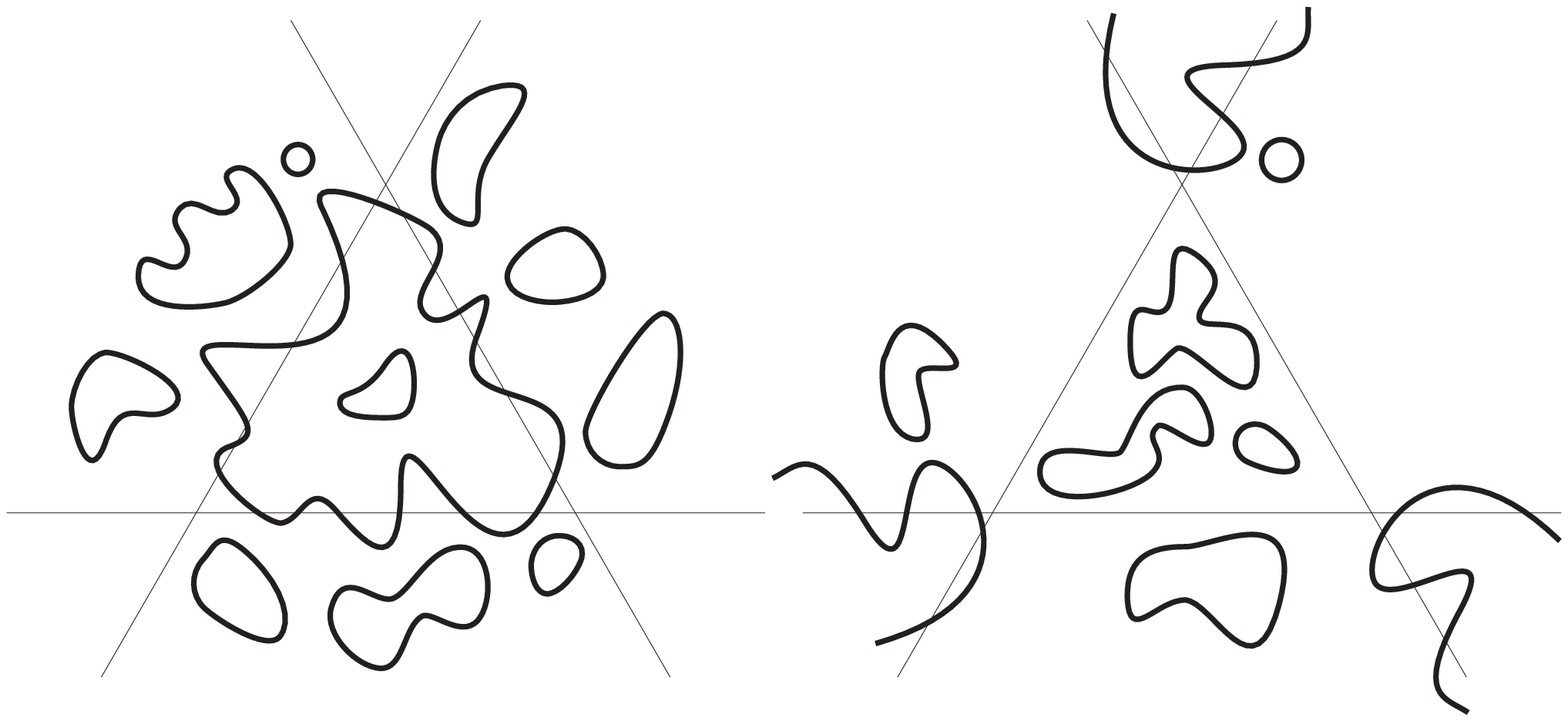}} 
    \caption{Oval configuration of  
    Harnack curves of degrees $6$ and $5$}
    \label{f2}
  \end{center}
\end{figure}

Among 
the $M$-curves, the Harnack curves are
distinguished by the position of their
ovals with respect to the each other and 
the  three coordinate
lines in $\P^2$ (they are also sometimes
called simple Harnack curves in the 
literature). The classical definition 
of a Harnack curve is somewhat complicated,
see \cite{Mi}. We will use instead modern
characterizations of Harnack curves 
obtained in \cite{Mi,MR}. 
By the main result of \cite{Mi}, for given 
degree $d$, a curve is a  Harnack curve
if and only if the topological configuration of ovals is as
illustrated in Figure \ref{f2}.
(It is different for odd and even degree.)
The numbers of ovals in each quadrant are 
consecutive triangular numbers. 
Of course, Figure \ref{f2} is
only meant to illustrate the topology and not 
be a plot of an actual Harnack curve. 

The group $(\R^\times)^2$ acts freely on the set of Harnack 
curves by rescaling the variables. We will call the 
quotient of this action the \emph{moduli space of 
Harnack curves}. We will see that the moduli space of 
degree $d$ Harnack curves is diffeomorphic to the 
closed octant $\R_{\ge0}^{(d+4)(d-1)/2}$, in particular, it is 
connected and simply connected. It follows that the
set of all Harnack curves has 4 connected components,
corresponding to a choice of one of 4 quadrants in 
Figure \ref{f2}. These choices are related by reversing 
the signs of $z$ or $w$. 

There is a corresponding 
flexibility on the dimer problem side: rescaling 
$z$ and $w$ corresponds to a natural operation
on dimer weights which was called \emph{changing
magnetic field} in \cite{KOS}. {}From many points
of view, it is natural to view this operation as
a generalized gauge transformation. For example, 
with such extended definition of gauge equivalence, 
any dimer model whose spectral curve is rational
is equivalent to an isoradial dimer, see Section \ref{sISO}.

\subsubsection{Amoeba map}\label{sH2}

Another useful characterization of Harnack 
curves $C$ is that the map
$$ 
(\C^*)^2\supset C(\C)\owns (z,w)\mapsto (\log|z|,\log |w|)
\in \cA(C) \subset \R^2
$$
from the curve $C$ to its amoeba $\cA(C)$ is
generically 2-to-1 over the interior of $\cA(C)$
(that is, 2-to-1 except at real nodes). 

The (geometric) genus $g$ of a Harnack curve is the number of 
compact ovals
that are not reduced to points. In 
particular, the amoeba of a genus $g$
Harnack curve has exactly $g$ compact 
holes, which is illustrated in Figure \ref{f5} 
for the case $g=3$. 

The 2-to-1 property 
constrains the 
possible singularities that $C$ can have
by constraining the possible links of the 
singular points. It was shown in \cite{MR}
that the only possible degenerations 
of Harnack curves with fixed Newton polygon
occur when some of the ovals shrink to 
zero size, producing real isolated double 
points (real nodes). In particular, 
it is impossible for 
two ovals of $C$ to meet.

\subsubsection{Ronkin function and Monge-Amp\`ere equation}

A third characterization of Harnack curves $C$,
obtained in \cite{MR} is by maximality of the
area of their amoeba $\cA(C)$. Concretely,
$C$ is Harnack if and only if its amoeba 
has the maximal possible area for given 
Newton polygon $\Delta$, namely $\pi^2$ 
times the area of $\Delta$. 

This maximality of the area is a consequence of 
the following Monge-Amp\`ere equation (see \cite{MR}) 
\begin{equation}
  \label{MA}
  \det 
  \begin{pmatrix}
    R_{xx} & R_{xy} \\
    R_{yx} & R_{yy} 
  \end{pmatrix}  = \frac1{\pi^2}\,,
\end{equation}
satisfied for all $(x,y)$ in 
the interior of the amoeba $\cA(C)$. Here 
$R(x,y)$ is the \emph{Ronkin function}  of 
$C$ defined by 
\begin{equation}
  \label{Ron}
  R(x,y) = \frac{1}{(2\pi i)^2} 
\iint_{\substack{
|z|=e^x\\
|w|=e^y}} \, 
\log |P(z,w)| \,\, \frac{dz}{z} \frac{dw}{w} \,. 
\end{equation}

\subsubsection{Holomorphic differentials}\label{speri}

Holomorphic differentials on $C$ will play 
an important role in our analysis. In this 
section we recall some basic facts about them. 

Suppose that $C$ has geometric genus $g$. 
Let $\alpha_1,\dots,\alpha_g$ be the 
compact real ovals of $C$ that are not reduced to points,
and let $\{q_i\}$,
$i=1,\dots,\binom{d-1}{2}-g$ denote its
isolated real nodes. The ovals $\alpha_i$ 
form a set of $a$-cycles in $H_1(C)$. The 
corresponding $b$-cycles can be taken in 
the form
\begin{equation}
  \label{bcycl}
  \beta_i =\{|z|=e^{x_i}, |w|<e^{y_i}\}\cap C
\end{equation}
where $(x_i,y_i)$ is a point inside the 
$i$th hole of the amoeba $\cA(C)$.
The
cycles $\beta_i$ are anti-invariant with the 
respect to the complex conjugation. 

Holomorphic differentials on $C$ have the 
form 
\begin{equation}
  \label{holodiff}
   \omega = \frac{Q(z,w)}
{\frac{\partial}{\partial w} P(z,w)} \, dz \,, 
\end{equation}
where $Q$ is a polynomial of degree $d-3$ 
vanishing at the points $\{q_i\}$, see \cite{ACGH}. It
is a fundamental fact about plane curves 
(see \cite{ACGH}, Appendix A) that the nodes 
$\{q_i\}$ impose independent conditions on 
polynomials of degree $d-3$. In particular, 
the space of differentials 
\eqref{holodiff} is $g$-dimensional. We 
choose a basis $\{\omega_i\}$ such that 
\begin{equation}\label{deltaij}
\int_{\alpha_i} \omega_j  = \delta_{ij} \,.
\end{equation}
The polynomials $Q_i$ are real and 
their $b$-periods are purely imaginary (since 
$\beta_i$ is anti-invariant under complex conjugation). 

The polynomial $Q_i$ has at least two zeros on
every compact 
oval except the $i$th, by the integral condition (\ref{deltaij}):
the form $dz/P_w$ is real and of constant sign on every oval, so 
$Q_i$ must change sign at least twice on each oval except the $i$th.
This implies that $Q_i$ must have no zeros on 
the $i$-th oval and precisely 2 zeros on any 
other oval, otherwise, the degree $d-3$ curve
$Q_i=0$ would intersect $C$ at 
$$
(d-1)(d-2) > d(d-3)
$$
points. It follows that the function 
\begin{equation}
  \label{intxy}
x \mapsto   \int_{x_0}^x \omega_i \,,
\end{equation}
where $x_0$ and $x$ are points on the
$i$th oval, is increasing and 
provides a natural parametrization of the 
$i$th oval. 

Let us call a degree $\binom{d-1}{2}$ divisor $D$ 
a \emph{standard divisor} if it has precisely one
point on each compact oval of $C$ (including degenerate ovals). 
By monotonicity of \eqref{intxy}, the product of 
compact ovals of $C$, which is the 
variety parameterizing  standard divisors, 
injects into the Jacobian $J(C)$
of degree $\binom{d-1}{2}$. It forms  
a component of its real locus.

\section{The spectral problem}

\subsection{Boundary of the spectral curve}\label{bdy}

By the \emph{boundary of a spectral curve} $C$
we mean the $3d$ points (counting multiplicity)
of its intersection with the coordinate lines
of $\P^2$. We show in this section that these points  have a 
simple interpretation in terms of the
weights of the dimer problem. In 
particular, they are real for real 
edge weights. 

\begin{figure}[hbtp]\psset{unit=0.5 cm}
  \begin{center}
    \begin{pspicture}(0,0)(10,6)
    \rput(5,3){\scalebox{1}{\includegraphics{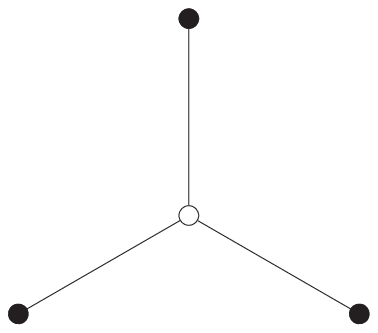}}} 
    \rput[rb](4.8,2.2){$v$}
    \rput[lB](7,1){$z a_v$}
    \rput[lB](5.1,4){$w b_v$}
    \rput[rB](3,1){$c_v$}
    \end{pspicture}
    \caption{The operator $\tbK(z,w)$}
    \label{f3}
  \end{center}
\end{figure}

It is useful to introduce the curve $\tC$
defined by the equation 
\begin{equation}
  \label{tP}
  \tP(z,w) = P(z^d,w^d) \,.
\end{equation}
This is a $d^2$-fold covering of the 
spectral curve $C$ ramified over its boundary. 
The advantage of going to this branched
covering is that operator $\bK(z^d,w^d)$
is gauge equivalent to the operator  $\tbK(z,w)$
defined by the rule illustrated in Figure \ref{f3}.
In the operator $\tbK(z,w)$,
the weights of \emph{all} vertical edges are multiplied by $w$ and weights of
all northwest-southeast edges are multiplied by $z$.

\begin{figure}[hbtp]\psset{unit=0.5 cm}
  \begin{center}
    \begin{pspicture}(-2,0)(12,10)
    \rput(5,5){\scalebox{0.64}{\includegraphics{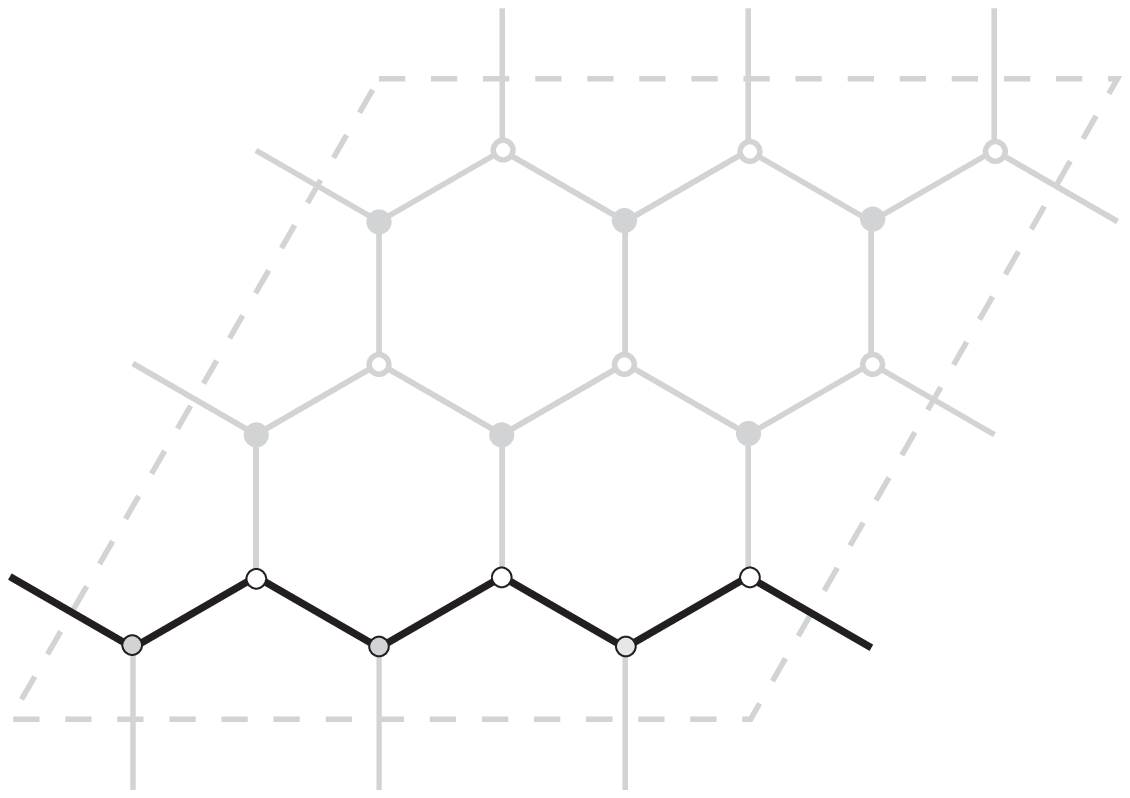}}} 
    \rput[lb](1.9,2.2){$a_1$}
    \rput[lb](5,2.2){$a_2$}
    \rput[lb](8.15,2.2){$a_3$}
    \rput[lt](0.1,2.2){$c_1$}
    \rput[lt](3.3,2.2){$c_2$}
    \rput[lt](6.5,2.2){$c_3$}
    \end{pspicture}
    \caption{The operator $\tbK(z,0)$}
    \label{f4}
  \end{center}
\end{figure}

In  particular, if $w=0$ then this operator
becomes block-diagonal with blocks corresponding
to chains of vertices of the form shown in 
Figure \ref{f4}.

The determinant of the corresponding block equals 
\begin{equation}
  \label{detch}
  \det 
  \begin{pmatrix}
    c_1  &  z a_2 && \\
& c_2 &  z a_3 \\
& & c_3 & \ddots  \\
 z a_d && & \ddots 
 \end{pmatrix} = \prod c_i + (-1)^{d+1} z^d \prod a_i \,.
\end{equation}

Therefore, the spectral curve $P(z,w)=0$ intersects the $w=0$ line 
at the points of the form
\begin{equation}
  \label{zbdry}
  z = (-1)^d \prod_{\textup{$v$ in a zig-zag cycle}} c_v/ a_v \,,
\end{equation}
where by a zig-zag cycle we mean a ``horizontal'' 
chain of vertices like in Figure \ref{f4}. 

The alternating products
of the from \eqref{zbdry} taken for all paths of all 
3 orientations multiply to $1$. This is, in fact,
a general constraint for the boundary points of any 
degree $d$ curve. It can be viewed, for example, as 
coming from the divisor class group of the union of
the coordinate lines. 

Note that once the fundamental 
domain is fixed, we get a canonical ordering of the 
boundary points of the spectral curve. It comes
from the ordering of the corresponding zig-zag cycles.

\subsection{The divisor of a vertex}

By definition, the spectral curve \eqref{spC} is where the 
kernel and the cokernel of the matrix $\bK(z,w)$ become nontrivial. 
For smooth spectral curves $C$, the kernel and the 
cokernel are $1$-dimensional \cite{CT}
and form a line bundle over $C$.  In particular, 
for every white vertex $v\in \Gamma$ 
the image of the function $\delta_v$ in 
$$
\Coker \bK(z,w)=\C^{\textup{white vertices}}/\Img \bK(z,w)
$$
defines a section of the 
cokernel bundle. We will denote by $(v)$ the 
divisor of this section excluding the points
on the boundary. 

In down to earth terms, 
equations of $(v)$ are the cofactors of 
the entries of $\bK(z,w)$ in the row 
corresponding to the white vertex $v$. Since
the columns of $\bK(z,w)$ are linearly dependent
on the curve $C$, $(v)$ is defined locally by 
the vanishing of a single minor. 

\begin{Theorem}\label{t1}
For real positive edge weights, the spectral 
curve $P(z,w)=0$ is 
a Harnack curve of degree $d$,  in particular, 
it has $\binom{d-1}{2}$  compact ovals. The divisor $(v)$
of any vertex $v$ is a standard divisor 
on this curve. 
\end{Theorem}

The first part of this theorem is one of the 
main results of \cite{KOS}. The argument given 
below provides a new proof of it.

\begin{proof}
It suffices to consider the case of generic edge weights,
since the properties of being a Harnack curve and
being a standard divisor are closed. 

Let us consider the bundle $\cL$ formed by 
the cokernels of the matrix $\tbK(z,w)$ on 
the curve $\tC$ defined by the equation \eqref{tP}. It fits into the 
following exact sequence 
\begin{equation}\label{exSeq}
  0 \to \cO(-1)^{d^2} \xrightarrow{\,\,\tbK\,\,} 
\cO^{d^2} \to \cL \to 0 
\end{equation}
of sheaves on $\P^2$. It follows that
$$
\chi(\cL)=d^2 \,.
$$
The curve $\tC$ has genus $(d^2-1)(d^2-2)/2$ and, therefore,
by Riemann-Roch we obtain 
$$
\deg \cL = \frac{d^2(d^2-1)}{2} \,.
$$
Let $v$ be a white vertex and consider the section $s_v$ of 
$\cL$ defined by the function $\delta_v$. The divisor 
of $s_v$ is the pull-back of $(v)$ except for 
the boundary contribution, which we will now determine. 

Let $p$ be a boundary point of $\tC$ corresponding to a
zig-zag cycle as in Figure \ref{f4} and suppose that the vertex $v$ is 
$k$ steps above this path. Below we prove the following: 

\begin{Lemma}\label{lem1} 
The order of the 
vanishing of $s_v$ at the point $p$ equals $k$. 
\end{Lemma}

Since every path 
gives rise to $d$ points on the boundary of $\tC$, 
it follows that the 
total contribution of the boundary points to the degree of $s_v$
equals $3d^2(d-1)/2$. Subtracting it, we obtain
\begin{equation}
  \label{deg(v)}
  \deg (v) = \frac1{d^2}(\deg\cL-\deg s_v)=\frac{(d-1)(d-2)}{2} \,.
\end{equation}

The remainder of the proof is based on the deformation to 
the constant weight case. If all edge weights are equal to 
one, the polynomial $\tP$ takes a particularly simple
form (see \cite{KOS})
\begin{equation}
  \label{tP0}
  \tP(z,w) = \prod_{i,j=1}^d (\varepsilon^i \,z + 
\varepsilon^j\, w + 1)\,,
\end{equation}
where $\varepsilon$ is a primitive $d$th root 
of unity. In particular, the spectral curve $C$
is a genus $0$ curve with $\binom{d-1}{2}$ 
isolated real nodes (isolated real solutions to
$\varepsilon^i \,z + 
\varepsilon^j\, w + 1=0$ exist when $i\neq j$ and $i,j\neq d$; moreover
replacing $i,j$ with $-i,-j$ gives an identical solution).
It also intersects each coordinate line once with 
multiplicity $d$, and, hence, satisfies the
topological definition of a Harnack curve.  

The space 
$$
\Ker \bK^*(z,w) = \left( \Coker \bK(z,w)\right)^*
$$
is $2$-dimensional
if $(z,w)$ is one of the nodes
of $C$. In particular, there exists an nonzero 
element of this space which annihilates $\delta_v$. 
By continuity, for all nearby curves
the divisor $(v)$ has a nearby
point, possibly complex. By \eqref{deg(v)} the total number of 
such points equals the degree of $(v)$, hence
all of them are simple real zeros (no point can bifurcate into
a complex conjugate pair of zeros).  

It follows that 
under a small real perturbation of weights, 
the spectral curve remains a Harnack curve and
$(v)$ remains a standard divisor. 
Because the ovals of a Harnack curve cannot meet
(see Section \ref{sH2}) and the boundary points 
do not become complex (see Section \ref{bdy}), the 
same statement holds globally for real positive 
weights.

\end{proof} 

\begin{proof}[Proof of Lemma \ref{lem1}] 
Let $(z_0,0)$ be the coordinates of the point $p$. 
Note that if the corresponding intersection of the 
curve $C$ with the $w=0$ axes was transverse (which 
happens generically) then $w$ is a local parameter 
on $\tC$ near $p$ and $z-z_0=O(w^d)$. Therefore, near $p$ we have
$$
\bK\big|_{\tC}=\begin{pmatrix}
  K_1 & w B_1\\
      & K_2 & w B_2\\
      &     & \ddots & \ddots\\
  w B_d    &     &     & K_d        
\end{pmatrix}+O(w^d)\,,
$$
where the blocks $K_i$ correspond to zig-zag paths as in \eqref{detch}
and $B_i$ are diagonal invertible matrices. The matrix $K_1$ has a 
one-dimensional kernel while all other $K_i$'s are invertible. 
Clearly,
$$
\Ker \bK^* = (v,-w K_2^{-1} B_1 v, w^2 K_3^{-1} B_2 K_2^{-1} B_1 v, \dots) + O(w^d)\,,
$$
where $v\in \Ker K_1^*$. Since the entries of $\Ker \bK^*$ are 
the cofactors of $\bK$ that we need, the lemma follows. 
\end{proof}

\subsection{Spectral transform}

We have constructed the following 
correspondence 
\begin{equation}
  \label{spTr}
  \left(
    \begin{gathered}
      \textup{edge weights}\big/\textup{gauge}\\
      \textup{fundamental domain} \\
      \textup{vertex} 
    \end{gathered} \right)
\quad \Rightarrow \quad 
\left(
\begin{gathered}
      \textup{Harnack curve}\\
      \textup{ordering of boundary points}\\
      \textup{standard divisor} 
    \end{gathered}
\right) \,,
\end{equation}
which we will call the \emph{spectral 
transform}. Our first result about it is
the following 

\begin{Theorem}\label{t2} 
  The spectral transform is injective. 
\end{Theorem}

\begin{proof}
  The spectral data determines the 
sheaf $\cL$ introduced in the proof 
of Theorem \ref{t1} up to isomorphism. 
By Theorem 1.1 in \cite{CT}, this 
determines the matrix $\tbK(z,w)$ up to 
multiplication on the left and on the 
right by elements of $GL(d^2,\R)$, that 
is, up to a choice of linear basis of 
both spaces in \eqref{Kas}.

Let $p$ be a boundary point of the spectral 
curve. It determines a filtration on the 
space $\C^{\textup{white vertices}}$ by 
the order of the vanishing of the corresponding
section of the cokernel bundle at $p$. It is 
easy to see that the basis $\{\delta_v\}$ is the 
unique, up to normalization, basis compatible
with all these filtrations. Same argument
applied to transpose matrix $\tbK^*$,
which has the same spectral curve, reconstructs
the delta-function basis of $\C^{\textup{black vertices}}$.
As a result we reconstruct  the matrix $\tbK(z,w)$ up to 
a multiplication on the left and on the right by 
a diagonal matrix, that is, up to a gauge
transformation. 
\end{proof}

A very useful (and well-known, cf.\ \cite{Be1,GNR,Hur,Skl}) 
way of thinking about the spectral data in \eqref{spTr}
is the following. The divisor 
is an ordered collection of distinct points of $\P^2$; 
similarly, the boundary points are an ordered collection 
of points on the coordinate axes. The spectral curve 
$P(z,w)$ is determined by this data as the unique degree $d$ curve 
passing through all these points. Indeed, since the 
boundary points are fixed, any other curve passing through
the same points is given by 
$$
P(z,w)+z w Q(z,w)=0\,, \quad \deg Q\le d-3\,.
$$ 
But then, as in Section \ref{speri}, $Q=0$ and 
$P=0$ will have $(d-1)(d-2)$ points in common, which 
is one too many. 

It is clear from this description that
the spectral data is parametrized by an open subset of 
$\R^{d^2+1}$ (which, in fact, is homeomorphic to $\R^{d^2+1}$, 
as we will see later). Also, this picture 
continues to work in a complex neighborhood 
and the injectivity remains valid with the same proof. 
It follows at once that: 

\begin{Corollary}
The differential of the spectral transform is injective. 
\end{Corollary}

\subsection{Surjectivity of the spectral transform}

This reconstruction procedure given in Theorem \ref{t2} 
can be made explicit to a certain degree, but it appears
to be  difficult to establish the 
surjectivity of the spectral transform 
directly. Instead, we will take a 
different route based on the following 

\begin{Proposition}
  The spectral transform is proper. 
\end{Proposition}

In fact, at several points in the paper we 
will resort to proving surjectivity by 
combining properness with injectivity of the 
differential. 

\begin{proof}
  Suppose that the spectral curve $C$ varies in some 
bounded set of Harnack curves. In particular, the 
the 
coefficients of $1$, $z^d$, and $w^d$ in the 
equation of $C$ are bounded from above and below in 
absolute value. 
Without loss of generality, we can assume that 
these coefficients are equal to $\pm 1$. This means 
that the weights of the 3 frozen matchings (matchings using
all edges of the same orientation) are 
equal to 1. Other 
coefficients of $P(z,w)$, which enumerate other 
periodic matchings, are bounded from above
in absolute value. This implies that the weight 
of any periodic matching is bounded from above. 

Periodic matchings are the vertices of the 
periodic flow polytope. This polytope is formed by nonnegative
periodic flows from white vertices to black vertices
such that the total flux is $1$ at each black vertex. The logarithm of the 
weight of a matching extends as a bounded linear
function to the interior of this polytope (recall that the 
weight of a matching is the product of its edge weights). 
The barycenter of this polytope is the 
flow of intensity $1/3$ along each edge. 
It has weight $1$ by our assumption. 

Coordinates on the space of weights 
modulo gauge are given by the combinations
of the form
\begin{equation}
  \label{loopweight}
  \prod_{i=1}^{k} \frac{\weight(e_{2i-1})}
{\weight(e_{2i})}
\end{equation}
where $e_1,e_2, \dots, e_{2k}$ 
is a closed loop of edges on the torus, 
such as for example the loop in Figure \ref{f4}. Such loops 
can be thought of as flows with zero flux. 
If the weight \eqref{loopweight}
of any such loop under the variation is unbounded, we can add a
small multiple of it to the constant $1/3$
flow and obtain a point inside the periodic flow 
polytope with unbounded weight. This 
contradiction completes the proof. 
\end{proof}

We will prove in Theorem \ref{d^2+1} that the spectral 
data on the left in \eqref{spTr} forms
a manifold diffeomorphic to $\R^{d^2+1}$. {}From the 
description \eqref{loopweight} of the 
coordinates on the space of weights 
modulo gauge, it it clear that this 
is also diffeomorphic to  $\R^{d^2+1}$.
Indeed, the cycle space  is
$d^2+1$-dimensional, generated by
the $d^2$ faces in a fundamental domain (subject to one relation),
and the horizontal and
vertical cycles. Therefore, we obtain 

\begin{Theorem}\label{t3}
The spectral transform is a bijection. 
\end{Theorem}

\section{The space of Harnack curves}
\label{s3}

\subsection{Harnack curves of genus zero}\label{Harnack}

Harnack curves of genus zero can be 
understood explicitly. By general theory, any real 
degree $d$ genus zero curve $C$ has a parametrization of the 
form 
\begin{equation}\label{Hg0}
z(t)=a_0\, \prod_{i=1}^d 
\frac{t-a_i}{t-c_i} \,, \quad 
w(t)=b_0\, \prod_{i=1}^d 
\frac{t-b_i}{t-c_i}\,,
\end{equation}
where $t\in \R P^1 = \R \cup \{\infty\}$, 
$a_0$ and $b_0$ are nonzero real numbers, 
and the other parameters are either real or appear in 
complex conjugate pairs. Such representation is 
unique up to the action of $PGL_2(\R)$
by fractional linear transformation of the 
variable $t$. 

The action of $PGL_2(\R)$ can by reduced down to 
the action of the connected group $PSL_2(\R)$ once we fix the orientation 
of the curve $C$. We will fix some cyclic ordering
of the coordinate axes and require that the  $C$ 
intersects them in that order (see Section \ref{sdefHar}). 
Let the cyclic order be: $\{z=0\}$, $\{w=0\}$, and the 
line at infinity.

\begin{Proposition}\label{pHg0}
The curve \eqref{Hg0} is Harnack if and only if the 
parameters in \eqref{Hg0} are real and cyclically 
ordered on $\R P^1$ as follows:
\begin{equation}
  \label{param1}
  a_1 \le a_2 \le \dots \le a_d < b_1 \le \dots \le b_d <
c_1 \le \dots \le c_d < a_1\,. 
\end{equation}
\end{Proposition}

\begin{proof}
In one direction, this is a part of the topological 
definition of the Harnack curve, see \cite{Mi}, which 
requires the curve $C$ to first intersect the line $\{z=0\}$ $d$ times, 
then the line $\{w=0\}$ also $d$ times, and, 
finally, the line at infinity $d$ times.

In the other direction, we have to show that as long as
the parameters \eqref{param1} are cyclically ordered 
the curve $C$ remains Harnack. Note that the strict inequalities
in \eqref{param1} ensure that $C$ doesn't pass
through the intersections of the coordinate axes, which in
other words means that the Newton polygon of its equation 
remains the full triangle.  Since
the real nodes of the curve cannot disappear or merge 
by the $2$-to-$1$ property (see Section \ref{sH2}), and the region in the 
parameter space defined by the inequalities \eqref{param1}
is connected, the result follows. 
\end{proof}

The $PSL_2(\R)$--quotient can be taken by fixing, for example,
\begin{equation}
  \label{SL2ga}
  a_1=0,\quad b_1=1, \quad c_1=\infty\,. 
\end{equation}
This gives the following: 

\begin{Corollary}
The moduli space of genus $0$ degree $d$
Harnack curves is diffeomorphic to $\R_{\ge 0}^{3d-3}$.  
\end{Corollary}

Let us associate to the curve $C$ its $3d$ points of 
intersection with the coordinate axes, counting 
multiplicity. For the curve \eqref{Hg0} these are
the points in projective coordinates $(0,A_i,1),(1,0,B_i),$ and $(C_i,1,0)$ where
\begin{equation}
  \label{ABC}
  A_i = b_0 \prod_j \frac{a_i-b_j}{a_i-c_j} \,,
\quad B_i = \frac 1{a_0} \prod_j \frac{b_i-c_j}{b_i-a_j}\,,
\quad C_i = \frac{a_0}{b_0} \prod_j \frac{c_i-a_j}{c_i-b_j}\,. 
\end{equation}
Note that all $A_i$ have the same sign (similarly for $B_i$'s and 
$C_i$'s) and the following relation: 
\begin{equation}
  \label{prodABC}
  \prod_i A_i \, B_i \, C_i = (-1)^d \,.
\end{equation}
Modulo the action of $(\R^\times)^2$, these constraints define a
manifold with boundary  diffeomorphic to $\R_{\ge 0}^{3d-3}$. To see this,
sort the $A_i$'s and look at the ratios between consecutive values. This 
$d-1$-tuple of '$A$' ratios is in $\R_{\ge 1}^{d-1}$. 
Similarly for the '$B$' and '$C$' ratios. 
For a given set of ratios the values of the smallest $A_i,B_i,C_i$ are then determined
up to $(\R^\times)^2$ by (\ref{prodABC}).
If on the other hand we choose an ordering of the 
boundary points as in \eqref{spTr} then the constraints define an $\R^{3d-3}$:
for example by the $(\R^\times)^2$-action we can fix $A_1=B_1=1$ then all
other $A_i,B_i,C_i$ except $C_1$ are free, $C_1$ being fixed by (\ref{prodABC}). 
 
We have the following:

\begin{Theorem}\label{tg0}
There exists a unique genus zero Harnack curve for
every choice of boundary points. 
\end{Theorem}

\begin{proof}
  The map \eqref{ABC} descends to a map 
\begin{equation}
  \label{abcABC}
    \{a_i,b_i,c_i\}_{i=1\dots d}\big/PSL_2(\R) 
\to \{A_i,B_i,C_i\}_{i=1\dots d}\big/(\R^\times)^2 \,,
\end{equation}
which becomes a map from $\R^{3d-3}$ to $\R^{3d-3}$ if we introduce
an ordering of points on both sides. It is easy to see that 
this map is proper. For example, in the gauge \eqref{SL2ga}, 
we clearly have $b_j-c_i > 1$ for any $j$ and any $i\ge 2$. 
If the ratio
$$
\frac{C_i}{C_1} = \prod_j \frac{a_j-c_i}{b_j-c_i}
$$
is bounded from below then 
$a_j-c_i$ is bounded from below. The same argument (using the fact that
$a_j-c_i$ and hence $-c_i$ is bounded below) shows that if 
$C_i/C_1$ and $A_i/A_1$ are bounded below then $b_j-a_i$ is bounded from below.
Finally if $B_i/B_1$ is bounded from below then
$$
\frac{B_i}{B_1} = \prod_{j=1}^d  \frac{1-a_j}{b_i-a_j}
\prod_{j=2}^d \frac{b_i-c_j}{1-c_j} \sim \frac1{b_i} \,, \quad  
b_i\gg 0 \,, 
$$
and we conclude that the $b_i$'s are bounded from above. 

It suffices to check, therefore, that the differential of the map 
\eqref{abcABC} is an isomorphism at every point. We compute
\begin{align}
\frac{\partial \log A_i}{\partial a_k} &= 
\delta_{ik} \sum_j \left(\frac1{a_i-b_j} + \frac1{c_j-a_i}\right) \notag\\
\frac{\partial \log A_i}{\partial b_j} &= - \frac1{a_i-b_j} = 
\frac{\partial \log B_j}{\partial a_i}\label{JlogA}
\,, \quad \textup{etc.}
\end{align}
Note that the Jacobi matrix \eqref{JlogA} is a symmetric matrix which 
is the sum of the following 
elementary $3\times 3$ blocks:  
\begin{equation}
  \label{block}
  \begin{pmatrix}
    \dfrac{1}{a_i-b_j} +  \dfrac{1}{c_k-a_i} &
-  \dfrac{1}{a_i-b_j}  &  - \dfrac{1}{c_k-a_i} \\
-  \dfrac{1}{a_i-b_j} &  \dfrac{1}{a_i-b_j} + \dfrac{1}{b_j-c_k}  &
-\dfrac{1}{b_j-c_k} \\
- \dfrac{1}{c_k-a_i} & -\dfrac{1}{b_j-c_k} & 
\dfrac{1}{c_k-a_i} + \dfrac{1}{b_j-c_k}
  \end{pmatrix}
\end{equation}
over all triples of indices
$i,j,k=1,\dots,d$. The matrix \eqref{block} has rank 1 
with kernel spanned by the vectors $(1,1,1)$ and $(a_i,b_j,c_k)$.
The remaining eigenvalue equals
$$
- \frac{(a_i-b_j)^2+(c_k-a_i)^2+(b_j-c_k)^2}{(a_i-b_j)(c_k-a_i)(b_j-c_k)}\,,
$$
which is nonvanishing and of 
the same sign for all triples $(i,j,k)$ by the cyclic
ordering condition \eqref{param1}. 

It follows that a vector is in the kernel of the Jacobi matrix
\eqref{JlogA} if and only if it is annihilated by every $3\times 3$
matrix \eqref{block}. It follows that this kernel is
spanned by $(1,1,\dots,1)$ and $$(a_1,\dots,a_d,b_1,\dots,b_d,c_1,\dots,c_d)$$ 
and coincides with the tangent space to the orbit 
of the subgroup $\mathfrak{B}\subset PSL_2(\R)$ formed by upper-triangular matrices. 
Indeed, replacing $t\mapsto\alpha t+\beta$ in (\ref{Hg0}) corresponds to the
change $a_i\mapsto (a_i-\beta)/\alpha$ and similarly for $b_i$ and $c_i$. 

Since the matrix \eqref{JlogA} is symmetric, its image is the 
orthogonal complement of its kernel. Now the tangent space to the 
$PSL_2(\R)$--orbit is $3$-dimensional, 
contains the kernel of the differential of \eqref{ABC}, 
and is mapped by this differential  to 
the tangent space to the $(\R^\times)^2$--orbit. 
We need to know the 
intersection of the image with the tangent space to the 
$(\R^\times)^2$--orbit. It is immediate to see that this 
intersection is always $1$-dimensional, in fact spanned by the vector
$$(\lambda,\dots,\lambda,\delta,\dots,\delta,-\lambda-\delta,\dots,-\lambda-\delta)$$
where $\lambda\sum a_i+\delta\sum b_i=(\lambda+\delta)\sum c_i.$
It follows that this one-dimensional space has to be the image of the tangent space to the 
$PSL_2(\R)$--orbit. 
This shows that the differential of 
\eqref{abcABC} is an isomorphism, which concludes the proof. 
\end{proof}

\subsection{Intercept coordinates}

Let $C$ be a Harnack curve with the equation $P(z,w)=0$. 
The Ronkin function \eqref{Ron} of the polynomial $P$
has a facet (that is, a domain on which is it affine 
linear) for every component of the amoeba complement, 
that is, for every monomial of the polynomial $P(z,w)$ except for 
degenerate components, where the facet is reduced to a point. 
The slope of the facet corresponding to the 
term $p_{ij}z^i w^j$ is always the same, namely $(i,j)$. 
The intercept, however, varies. We will now show that
these intercepts can be taken as local 
coordinates on the space of Harnack curves. 

Since the 
intercepts of the unbounded components are easily 
found from the coefficients on the boundary of the
Newton polygon of $P$, and hence, from the boundary 
points of $C$, we can alternatively take 
the points on the boundary and the intercepts of the 
bounded facets as our coordinates. 

The space of Harnack curves with given boundary data and
genus at most $g$ is a real semialgebraic variety of dimension 
$g$, naturally stratified by the possible genus
degenerations (where one or more ovals 
degenerates to a point). We have the following 

\begin{Proposition}\label{prop4}
  The variety of Harnack curves with given boundary 
and genus is smooth with local coordinates 
given by the intercepts of the nontrivial compact ovals. 
\end{Proposition}

\begin{proof}
The tangent space to a curve $C$ with given boundary and nodes
is formed by the space of 
polynomials $R$ of degree $d$ vanishing at all nodes
and boundary points of $C$, modulo the polynomial $P$ itself. 
Such polynomials have the form $zw \, Q(z,w)$ where $Q$
is a polynomial of degree $\le d-3$ vanishing at the nodes of $C$.
It follows (see Section \ref{speri})
that the space of such polynomials is $g$-dimensional. 

Let $Q$ be such polynomial and let $(x,y)$ be a point inside the $k$th hole of 
the amoeba. We compute the variation of the $R(x,y)$ in the 
direction $Q$ as follows 
\begin{multline}
  \frac{d}{dt}\, R_{P+tzw\,Q}(x,y)\Big|_{t=0}\, 
 =  \frac{1}{(2\pi i)^2} 
\iint_{\substack{
|z|=e^x\\
|w|=e^y}} \, 
\frac{Q(z,w)}{P(z,w)} \,\, dz \, dw 
=\\
\frac{1}{2\pi i}  \int_{|z|=e^x}  dz 
\sum_{\substack{
P(z,w_r)=0\\
|w_r|< e^y}} \frac{Q(z,w_r)}
{\frac{\partial}{\partial w}P(z,w_r)} 
= \frac{1}{2\pi i}  \int_{\beta_k} \frac{Q(z,w)}
{\frac{\partial}{\partial w}P(z,w)} \, dz \,,
\end{multline}
where $R_{P+tzw\,Q}$ denotes the Ronkin function of the 
polynomial $P+tzw\,Q$, the integral over $w$ is computed
by residues, and $\beta_k$ is the $b$-cycle corresponding to the 
$k$th hole, see Section \ref{speri}. 

We see that the Jacobian of the transformation from 
the coefficients to the intercepts is the period matrix of 
the curve $C$, and, in particular, it is  nondegenerate. 
\end{proof}

Proposition \ref{prop4} can be strengthened as follows:

\begin{Proposition}
The neighborhood of genus $g$ degree $d$ Harnack curve inside the space 
of all  degree $d$ Harnack curves with the same boundary is diffeomorphic to 
$$
\R^g \times (\R_{\ge 0})^{\binom {d-1}2-g}\,,
$$
where the genus $g$ stratum is embedded as $(\R^g,0)$. 
\end{Proposition}

\begin{proof}
The $\R_{\ge 0}$ coordinates are given by 
the values of the polynomial $Q$ from the proof of Proposition \ref{prop4} 
at the nodes of $C$. These values have to be of definite sign 
in order for an oval to be present but otherwise arbitrary by
Section \ref{speri}.
\end{proof}

\subsection{Variational principle}

The holes of the amoeba satisfy the following 
variational principle which is parallel to 
the variational principle for conformal maps. 

\begin{Proposition}\label{pva}
  Under a boundary preserving variation that lowers one intercept
while keeping all others fixed, the corresponding
hole in the amoeba shrinks while all other 
components of the amoeba complement, including the 
unbounded ones, expand. 
\end{Proposition}

The proof will use the Legendre transform 
$\Rc$ of the Ronkin function. In the dimer
problem, it has the meaning of the 
\emph{surface tension} function. It is defined 
on the Newton polygon of $P$,  satisfies the 
Monge-Amp\`ere equation
\begin{equation}
  \label{MAc}
  \det 
  \begin{pmatrix}
    \Rc_{xx} & \Rc_{xy} \\
    \Rc_{yx} & \Rc_{yy} 
  \end{pmatrix}  = \pi^2\,,
\end{equation}
and has conical singularities (which correspond to facets of $R$)
at the lattice points inside the Newton polygon corresponding to
complementary components . 

\begin{proof}
  Let $\delta\Rc$ denote the variation of $\Rc$. It is 
positive at one singularity of $\Rc$ and zero 
at all other singularities. It also vanishes on 
the boundary of the Newton polygon. 

The function $\delta\Rc$
satisfies the equation 
$$R_{xx}\delta R_{yy}+R_{yy}\delta R_{xx}-R_{xy}\delta R_{yx}-R_{yx}\delta R_{xy}=0,$$
which is the linearization of \eqref{MAc}, and can be written
\begin{equation}
  \label{linMAc}
  \tr \,  (\partial^2 \Rc)^{-1} (\partial^2 \delta \Rc) = 0.
\end{equation}
 Here 
$\partial^2$ denotes the Hessian matrix of a function. 
Away from singularities, the function $\Rc$ is strictly 
convex, so $\partial^2\Rc$ is positive definite,
 and hence \eqref{linMAc} satisfies a maximum 
principle: local extrema of $\delta \Rc$ are not allowed. 
{}From the boundary conditions, we conclude that $\delta \Rc\ge 0$
everywhere. It follows that the cones at all vertices 
except one are becoming more acute, hence all but one of the
facets of $R$ are expanding. 

Adding a constant to $\Rc$, which doesn't affect \eqref{MAc},
we can achieve that $\delta\Rc$ is zero at one singularity and 
negative at all others. It follows that the cone at 
this one singularity is becoming more obtuse and, hence, 
the corresponding facet of $R$ is shrinking. 
\end{proof}

The following Figure \ref{f5}
 (best seen in color) illustrates the variational principle: 
\begin{figure}[!hbtp]
  \begin{center}
    \scalebox{0.4}{\includegraphics{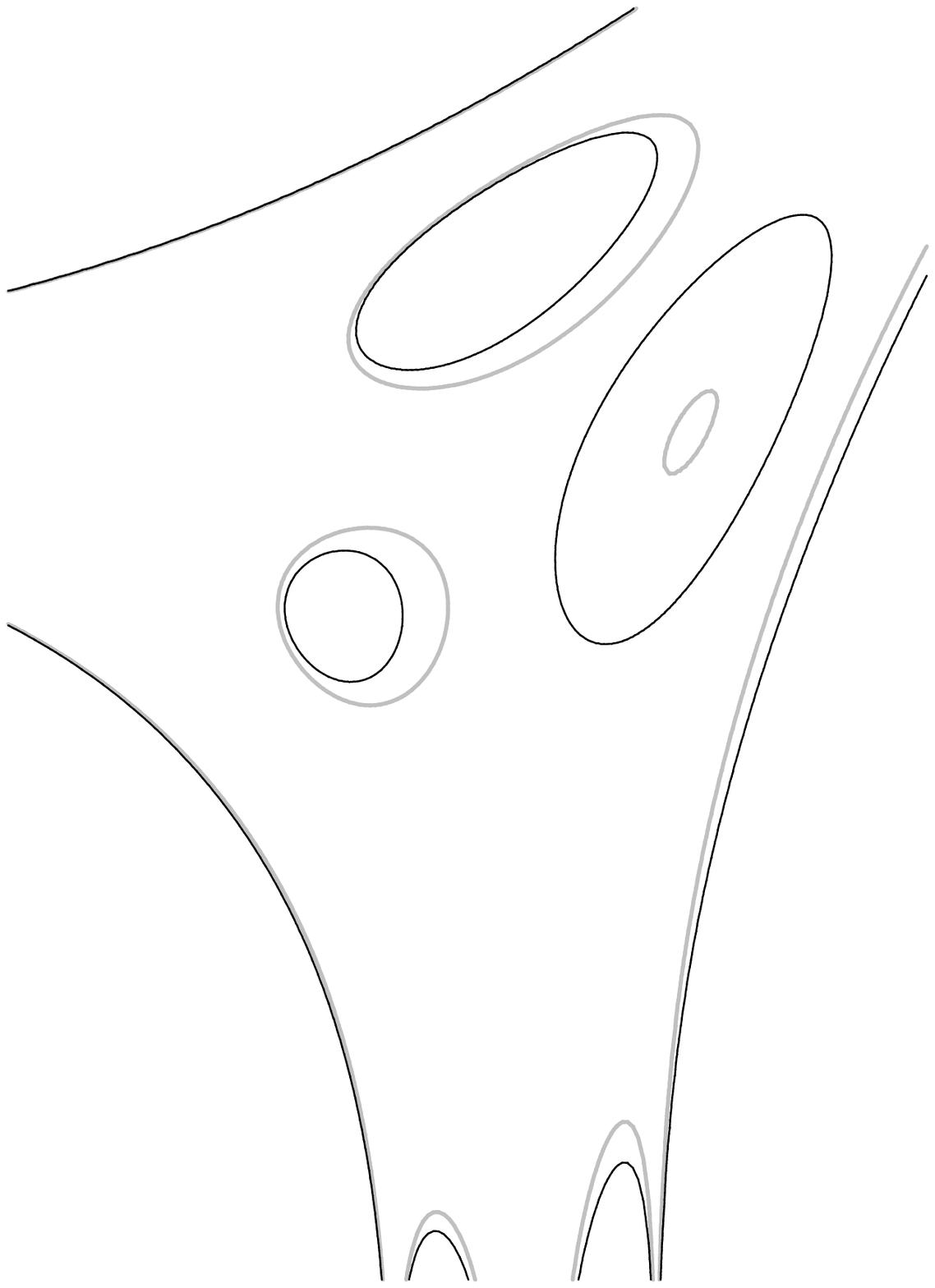}} 
    \caption{One oval is shrinking, the others expanding.}
    \label{f5}
  \end{center}
\end{figure}

\subsection{Volume under the Ronkin function}

The maximum principle for the linearized Monge-Amp\`ere equation
used in the proof of Proposition \ref{pva} has another application
for the volume under the Ronkin function.

Let $C_1$ and $C_2$ be two Harnack curves with the same boundary
data 
and assume that their equations $P_1$ and $P_2$ are normalized
in the same way (such as e.g. both constant terms are equal to 
1). Then the corresponding Ronkin functions $R_1$ and $R_2$  
agree asymptotically at infinity and the integral 
$$
\Vol R_1 - \Vol R_2 = \iint_{\R^2} (R_1 - R_2) \, dx dy
$$
converges. Choosing any $R_2$ as our reference point, this 
defines a functional $\Vol R$ of a Harnack curve. This function 
has the following monotonicity in the 
intercept coordinates. 

\begin{Proposition}\label{volmin}
  A variation of the curve $C$ which lowers the 
intercepts also lowers $\Vol C$. The genus zero Harnack 
curve is the unique volume
minimizer with given boundary conditions. 
\end{Proposition}

\begin{proof} This follows from the same maximum principle 
used in the proof of Proposition \ref{pva}: a variation which lowers the
intercepts necessarily lowers the entire amoeba. The intercepts can be lowered
to the unique (by Theorem \ref{tg0}) 
point where the amoeba holes shrink to points, that is,
to the genus-zero curve.
\end{proof}

\begin{Theorem}\label{d^2+1}
  The variety of spectral data in \eqref{spTr} is 
diffeomorphic to $\R^{d^2+1}$. 
\end{Theorem}

\begin{proof}
The gradient flow of the volume functional with 
respect to any metric on the projective space of 
degree $d$ curves contracts the variety of Harnack 
curve with given boundary to a small neighborhood 
of the genus zero curve. It follows that this 
space, together with its stratification by the 
genus is diffeomorphic to the product
of $\binom{d-1}{2}$ copies of $\R_{\ge 0}$. Adding 
a standard divisor makes it diffeomorphic  
to $\R^{(d-1)(d-2)}$ (each hole, possibly degenerate,
with a point on its
boundary is diffeomorphic to an $\R^2$). Adding in the boundary
points, which we showed in section
\ref{Harnack} was a space diffeomorphic to $\R^{3d-1}$,
makes $\R^{d^2+1}$.
\end{proof}

\subsection{Areas of holes}

Under the variation considered in Proposition \ref{pva}
the area of one of the holes in the amoeba is decreasing,
while the area of all other ones is increasing. Since 
the outer components of the amoeba complement are also 
expanding and the total area of the amoeba is 
preserved, we conclude that the decrease in the area of
one hole dominates the total increase in the areas of 
all other holes. In other words, the Jacobian of 
the transformation from intercepts to the areas of 
holes in the amoeba is strictly diagonally dominant and, 
hence, nonsingular. We obtain the following

\begin{Proposition}\label{parhlc}
  The areas on the amoeba holes can be takes as local coordinates
on the manifold of Harnack curves with given boundary and 
genus. 
\end{Proposition}

This can be sharpened as follows: 

\begin{Theorem}
The areas of the amoeba holes map the set of all 
Harnack curves with given boundary diffeomorphically 
onto $\R_{\ge 0}^{(d-1)(d-2)/2}$, mapping the 
stratification by the genus to the stratification 
by the number of nonzero coordinates. 
\end{Theorem}

\begin{proof}
It remains to show that the map from the curve $C$
with fixed boundary to the areas of its amoeba holes is proper.
Suppose that some coefficients of the equation 
$P(z,w)=\sum_{i,j} p_{ij}\, z^i w^j$ of $C$ 
grow to infinity. Recall that fixing boundary means
fixing the coefficients on the boundary $\partial\Delta$ 
of the Newton polygon $\Delta$, so only the interior coefficients
$p_{ij}$ can grow. 

For each point $(i,j)$ in the interior of $\Delta$ define the
convex set
$$C_{i,j}=\{(a,b)\in\R^2~|~
ai+bj+\log p_{ij} \ge  \max_{(k,l)\ne (i,j)}\big(ak+bl+\log p_{kl}\big)\}.$$
The union over $(i,j)$ of the $C_{i,j}$ is the set
\begin{equation}
  \label{dualch}
 \{(a,b)\in\R^2~|~ \max_{(i,j)\in\Delta}\big(ai+bj+\log p_{ij}\big) \ge
\max_{(i,j)\in \partial \Delta}\big(ai+bj+\log p_{ij}\big)\}.
\end{equation}
By our assumption, this set grows to infinity, in 
the sense that it eventually contains a ball of any arbitrarily large 
radius. Therefore
one of the sets $C_{i,j}$ eventually contains a ball of 
radius $r$, where $r$ can be arbitrarily large. 

Since the Euclidean distance between $(i,j)$ and any lattice point 
$(k,l)\ne (i,j)$ is at least $1$, for any $(k,l)\ne (i,j)$
there exists $(a_0,b_0)$ satisfying $a_0^2+b_0^2\le (r/2)^2$ such that
$$
a_0(k-i)+b_0(l-j) \ge r/2.
$$ 
It follows that if $C_{i,j}$ contains a ball of radius $r$,
the set of $(a,b)$ satisfying
\begin{equation}
 \label{dualpa2}
 ai+bj+\log p_{ij} >  r/2 + \max_{(k,l)\ne (i,j)}\big(ak+bl+\log p_{kl}\big)\,,
\end{equation}
contains a ball of radius $r/2$. The inequality \eqref{dualpa2} 
implies that 
$$
|p_{ij}\,z^i w^j| > e^{r/2} \max_{(k,l)\ne (i,j)} |p_{kl}\,z^k w^l|\,,
\quad |z|=e^a\,, |w|=e^b\,,
$$
therefore, all points $(a,b)$ satisfying \eqref{dualpa2} lie in the $(i,j)$ 
component of the amoeba complement provided $r$ is large enough. 
It follows that the area of that component is unbounded, which 
concludes the proof. 
\end{proof}

\begin{Corollary}
The areas of the amoeba holes and the distances between the amoeba
tentacles provide global coordinates on the moduli space of 
Harnack curves.   
\end{Corollary}

Note that the distances between the amoeba
tentacles can be viewed as ``renormalized'' 
areas of the semibounded components of the 
amoeba complement.





\section{Isoradial dimers and genus zero curves}\label{sISO}

Since Harnack curves of genus zero played a 
special role in our analysis, it is natural 
to ask for a characterization of those dimer
weights that lead to spectral curves of 
genus zero. Genus zero weights are distinguished, 
for example, by maximizing the partition 
function for given boundary conditions. 

In this section, we show that, up to gauge
equivalence, genus zero weights are the same
as isoradial weights studied in \cite{Ke2}. 
The results of \cite{Ke2} yield an explicit
rational parametrization \eqref{parsin} of the spectral 
curve for isoradial dimers.
In genus zero, the spectral curve determines
the weights uniquely, up to gauge transformation. 
It remains, therefore, to characterize the curves
\eqref{parsin} inside all genus zero Harnack 
curves, which is the content of Proposition \ref{isorad}

\subsection{Isoradial embeddings}

\begin{figure}[!hbtp]
  \begin{center}
    \scalebox{0.4}{\includegraphics{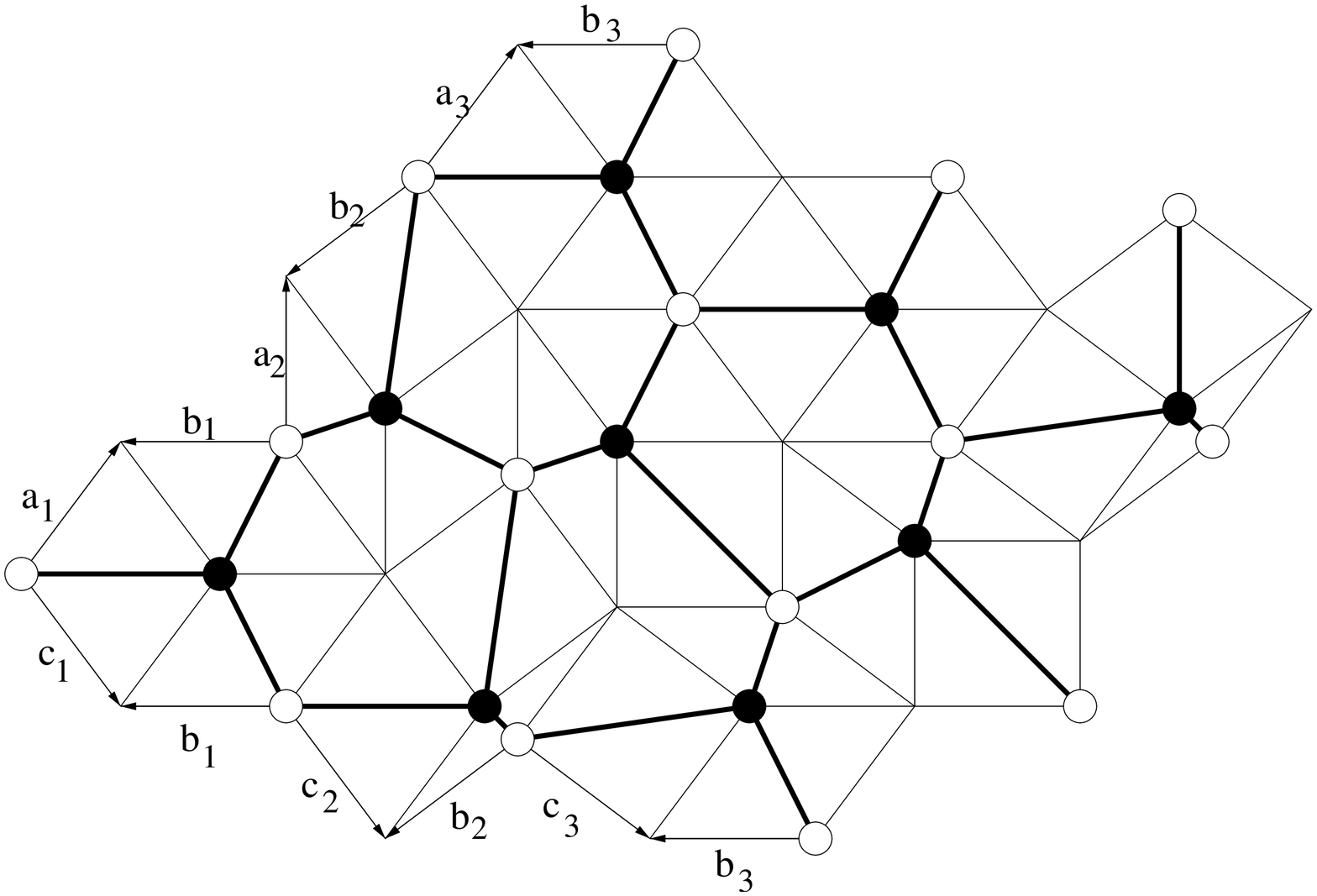}} 
    \caption{Isoradial embedding of the honeycomb graph.}
        \label{isoradhex}
  \end{center}
\end{figure}

If $\Gamma$ is embedded so that every face
is inscribed in a circle of radius $1$,
we say that it is {\it isoradial} on condition that
the weight of an edge is $\sqrt{4-\ell^2}$,
where $\ell$ is its length. Thus the weight is the 
distance between the circumcenters of the two 
faces adjacent to the edge. See Figure \ref{isoradhex}.

For the hexagonal lattice, an isoradial embedding
with $d\times d$ fundamental domain is determined by three 
$d$-tuples of unit modulus complex numbers,
$\{a_1,\dots,a_d\},\{b_1,\dots,b_d\}$ and $\{c_1,\dots,c_d\}$,
which are edges of the rhombi connecting vertices to
the face centers as indicated in 
Figure \ref{isoradhex}. They obey the condition that
they lie in three disjoint subintervals of the circle, that is,
going counterclockwise
around the circle one encounters first the $a$'s then the $b$'s then the $c$'s.
It is natural to consider these parameters up to 
a simultaneous rotation. 

\subsection{Bloch-Floquet eigenfunctions}

In \cite{Ke2} it was shown that
Bloch-Floquet eigenfunctions on isoradial graphs
satisfy the ``local multiplier condition'',
which says that
if $\b$ and $\w$ are (respectively black and white) adjacent vertices on an edge
bounding faces with circumcenters $C_1$ and $C_2$,
then
$$f(\b)=\frac{f(\w)\sqrt{r_1r_2}}{(u-r_1)(u-r_2)}$$
where $u$ is a complex parameter and $r_i=C_i-\w$ are unit complex
numbers.  The constant $\sqrt{r_1r_2}$ does not appear in \cite{Ke2}
due to the fact that we are using a different gauge for the Kasteleyn matrix.

In particular the eigenfunction has Floquet multipliers
\begin{eqnarray}\label{mults}
z(u) &=& \frac{f(\w+(1,0))}{f(\w)}=\prod_{i=1}^d\frac{u-a_i}{u-b_i}\sqrt{\frac{b_i}{a_i}}\\
w(u)&=&\frac{f(\w+(0,1))}{f(\w)}=\prod_{i=1}^d\frac{u-c_i}{u-b_i}\sqrt{\frac{b_i}{c_i}}. 
\notag 
\end{eqnarray}
The parameter $u\in\P^1$ parametrizes the spectral curve $P(z,w)=0$
which is, therefore, rational. Setting
$$
u=e^{it}\,, \quad a_i=e^{i\alpha_i}\,, \quad b_i=e^{i\beta_i}\,, \quad 
c_i=e^{i\gamma_i}\,,
$$
the parametrization \eqref{mults} becomes the following 
\begin{equation}
  \label{parsin}
  z(t) = \prod_i \frac{\sin\frac{t-\alpha_i}{2}}{\sin\frac{t-\beta_i}{2}}\,,
\quad
w(t) = \prod_i \frac{\sin\frac{t-\gamma_i}{2}}{\sin\frac{t-\beta_i}{2}}\,,
\end{equation}
in particular, $t\in \R$ parametrizes the unique nontrivial oval.

\subsection{Characterization of isoradial curves}

Compare \eqref{parsin} to \eqref{Hg0} and note the absence of arbitrary constant 
factors in front. We will call curves of the form \eqref{parsin} the 
\emph{isoradial curves}. They have the following simple characterization 

\begin{Proposition}\label{isorad}
  A genus zero Harnack curve $C$ is isoradial if and only if the origin is
in the amoeba of $C$. 
\end{Proposition}

\begin{proof}
Putting $u=0$ in  \eqref{mults} we obtain a point on the unit torus, 
so the origin is in the amoeba of any such curve. 

Conversely, making the standard change of variable from the 
upper half plane to unit disk in Proposition \ref{pHg0}, we see that 
any Harnack curve of genus zero will have a parametrization of almost
the same form as \eqref{mults}, except we have to allow two arbitrary real 
factors in front. Such parametrization is unique up to the 
$SL_2(\R)$--action by automorphism of the unit disk, so we need to check whether
the $SL_2(\R)$--action can be used to set these arbitrary factors to $1$. 
Conversely, we can act by $SL_2(\R)$ on \eqref{mults} and see what 
kind of constant factors we can get. 

Since rotating the unit disk around the origin clearly does not change 
anything, we have to look at the transformation 
$$
T(u) = \frac{u-\zeta}{1-\bar \zeta u}\,, \quad |\zeta|< 1 \,.
$$
Applying it to \eqref{mults}, we obtain the parametrization: 
\begin{align*}
  \tilde z(u) & =\prod_{i=1}^d \frac{|\zeta+a_i|}{|\zeta+b_i|} \, 
\prod_{i=1}^d\frac{u-\tilde a_i}{u-\tilde b_i}
\sqrt{\frac{\tilde b_i}{\tilde a_i}}\\
\tilde w(u)&= \prod_{i=1}^d \frac{|\zeta+c_i|}{|\zeta+b_i|} \, 
\prod_{i=1}^d\frac{u-\tilde c_i}{u-\tilde b_i}\sqrt{\frac{\tilde b_i}{\tilde c_i}}\,,
\end{align*}
where $\tilde a_i= T^{-1}(a_i)$ and similarly for $\tilde b_i$ and $\tilde c_i$. 
Comparing this with \eqref{mults}, we observe that
$$
\prod_{i=1}^d \frac{|\zeta+a_i|}{|\zeta+b_i|} = |z(-\zeta)|\,, \quad
\prod_{i=1}^d \frac{|\zeta+c_i|}{|\zeta+b_i|} = |w(-\zeta)|\,, 
$$
and hence, by its very definition, the amoeba of the curve $C$ describes the 
possible values of these factors. 
\end{proof}

Note from the proof that the $2$-to-$1$ property implies that 
point $\zeta$, when it exists, is uniquely defined. 

By a simple rescaling of $z$ and $w$ we can shift the 
amoeba of any curve $C$ so that it contains the zero. In particular, 
for any genus-$0$ Harnack curve this gives a canonical 
family of isoradial parameterizations indexed by the 
points in the amoeba.

\end{document}